\documentclass[10pt]{article}
\usepackage{authblk}
\usepackage{graphicx}
\usepackage{epstopdf}

\newtheorem{theorem}{Theorem}

\newtheorem{corollary}[theorem]{Corollary}

\begin{document}

\title{The algebraic properties of the space- and time-dependent one-factor
model of commodities}
\author[1]{A Paliathanasis\thanks{%
paliathanasis@na.infn.it}}
\author[2]{RM Morris\thanks{%
rmcalc85@gmail.com}}
\author[2,3,4]{ PGL Leach\thanks{%
leach@ucy.ac.cy}}

\affil[1]{Instituto de Ciencias F\'{\i}sicas y Matem\'{a}ticas, Universidad Austral de Chile, Valdivia, Chile}
\affil[2]{Department of Mathematics and
Institute of Systems Science, Research and Postgraduate Support, Durban
University of Technology, PO Box 1334, Durban 4000, Republic of South Africa}
\affil[3]{School of Mathematics, Statistics and Computer Science, University
of KwaZulu-Natal, Private Bag X54001, Durban 4000, Republic of South Africa}
\affil[4]{Department of Mathematics and Statistics, University of Cyprus,
Lefkosia 1678, Cyprus}

\renewcommand\Authands{ and }
\maketitle

\begin{abstract}
We consider the one-factor model of commodities for which the parameters of
the model depend upon the stock price or on the time. For that model we
study the existence of group-invariant transformations. When the parameters
are constant, the one-factor model is maximally symmetric. That also holds
for the time-dependent problem. However, in the case for which the
parameters depend upon the stock price (space) the one-factor model looses
the group invariants. For specific functional forms of the parameters the
model admits other possible Lie algebras. In each case we determine the
conditions which the parameters should satisfy in order for the equation to
admit Lie point symmetries. Some applications are given and we show which
should be the precise relation amongst the parameters of the model in order
for the equation to be maximally symmetric. Finally we discuss some
modifications of the initial conditions in the case of the space-dependent
model. We do that by using geometric techniques.
\end{abstract}


\section{Introduction}

Three models which study the stochastic behaviour of the prices of
commodities that take into account several aspects of possible influences on
the prices were proposed by E Schwartz \cite{schwartz} in the late nineties.
In the simplest model (the so-called one-factor model) Schwartz assumed that
the logarithm of the spot price followed a mean-reversion process of
Ornstein--Uhlenbeck type. \ The one-factor model is expressed by the
following $(1+1)$ evolution equation%
\begin{equation}
\frac{1}{2}\sigma ^{2}S^{2}F_{,SS}+\kappa \left( \mu -\lambda -\log S\right)
SF_{,S}-F_{,t}=0,  \label{1fm.01}
\end{equation}%
where $\kappa >0$ measures the degree of mean reversion to the long-run mean
log price, $\lambda $ is the market price of risk, $\sigma $ is the standard
deviation of the return on the stock, $S$ is the stock price, $\mu $ is the
drift rate of $S$ and $t$ is the time. $F$ is the current value of the
futures contract which depends upon the parameters $t,~S,$, \textit{i.e.}, $%
F=F\left( t,S\right) $.

Generally $\kappa $, $\lambda $, $\sigma $ and $\mu $ are assumed to be
constants. In such a case the closed-form solution of equation (\ref{1fm.01}%
) which satisfies the initial condition%
\begin{equation}
F\left( 0,S\right) =S  \label{1fm.01a}
\end{equation}%
was given in \cite{schwartz}. It is
\begin{equation}
\ln F\left( t,S\right) =e^{-\kappa t}\ln S+\left( 1-e^{-\kappa t}\right)
a^{\ast }+\frac{\sigma ^{2}}{4\kappa }\left( 1-e^{-2\kappa t}\right)
\label{1fm.02}
\end{equation}%
with $a^{\ast }=\mu -\lambda -\frac{1}{2}\frac{\sigma ^{2}}{\kappa }$.

It has been shown that the closed-form solution (\ref{1fm.02}) follows from
the application of Lie point symmetries. In particular it has been shown
that equation (\ref{1fm.01}) is of maximal symmetry, which means that it is
invariant under the same group of invariance transformations (of dimension $%
5+1+\infty $) as that of the Black-Scholes and the Heat Conduction Equation
\cite{IbragB}. The detailed analysis for the Lie symmetries of the three
models, which were proposed by Schwartz, and the generalisation to the $n$%
-factor model can be found in \cite{Leach08}. Other Financial models which
have been studied with the use of group invariants can be found in \cite%
{Leach05a,Leach06a,Naicker,Sinkala08a,Sinkala08b,wafo,consta,Lescot,dimas2}
and references therein.

Solution (\ref{1fm.02}) is that which arises from the application of the
invariant functions of the Lie symmetry vector
\begin{equation}
X_{sol}= \mbox{\rm e}^{\kappa t}S\partial _{S}+F\partial _{F}  \label{1fm.03}
\end{equation}%
and also leaves the initial condition invariant.

In a realistic World parameters are not constants, but vary in time and
depend upon the stock price, that is, the parameters have time and
\textquotedblleft space\textquotedblright\ dependence \cite{finBook1,finBook2}, where as space we
mean the stock price parameters as an analogue to Physics. In this work we
are interested in the case for which the parameters $\kappa $, $\lambda $, $%
\sigma $ and $\mu $ are space dependent, \textit{ie}, are functions of $S$.
We study the Lie point symmetries of the space-dependent equation (\ref%
{1fm.01}). As we see in that case, when $\sigma _{,S}\neq 0$, there does not
exist any Lie point symmetry which satisfies the initial condition (\ref%
{1fm.01a}). The Lie symmetry analysis of the time-dependent
Black-Scholes-Merton equations was carried out recently in \cite{time}, it
has been shown that the autonomous, and the nonautonomous
Black-Scholes-Merton equation are invariant under the same group of
invariant transformations, and they are maximal symmetric. The plan of the paper is as follows.

The Lie point symmetries of differential equations are presented in Section %
\ref{Preliminaries}. In addition we prove a theorem which relates the Lie
point symmetries of space-dependent linear evolution equations with the
Homothetic Algebra of the underlying space which defines the Laplace
operator. In Section \ref{space1} we use these results in order to study the
existence of Lie symmetries of for the space-dependent one-factor model (\ref%
{1fm.01}) and we show that the space-dependent problem is not necessarily
maximally symmetric. The generic symmetry vector and the constraint
conditions are given and we prove a corollary in with the space-dependent
linear evolution equation is always maximally symmetric when we demand that
there exist at least one symmetry of the form (\ref{1fm.03}) which satisfies
the Schwartz condition (\ref{1fm.01a}). Furthermore in section \ref{proof2}
we consider the time-dependence problem and we show that the model is always
maximally symmetric. Finally in section \ref{con} we discuss our results and
we draw our conclusions. Appendix\ref{proof1} completes our analysis.

\section{Preliminaries}

\label{Preliminaries}

Below we give the basic definitions and properties of Lie point symmetries
for differential equations and also two theorems for linear evolution
equations.

\subsection{Point symmetries of differential equations}

By definition a Lie point symmetry,~$X,$ of a differential equation \newline
\[
\Theta\left( x^{k},u,u_{,i},u_{,ij}\right) =0,
\]
where the $x^{k}$ are the independent variables, $u=u\left( x^{k}\right) $
is the dependent variable and
\[
u_{,i}=\frac{\partial u}{\partial x^{i}}
\]
is the generator of a one-parameter point transformation under which the
differential equation $\Theta $ is invariant.

Let $\left( x^{i},u\right) \rightarrow \left( \bar{x}^{i}\left(
x^{k},u,\varepsilon \right) ,\bar{u}\left( x^{k},u,\varepsilon \right)
\right) $ be a one-parameter point transformation of the independent and
dependent variables with the generator of infinitesimal transformations
being
\begin{equation}
X=\xi ^{i}\left( x^{k},u\right) \partial _{i}+\eta \left( x^{k},u\right)
\partial _{u}.  \label{go.10}
\end{equation}

The differential equation $\Theta $ can be seen as a geometric object on the
jet space $J=J\left( x^{k},u,u_{,i},u_{,ij}\right) $. Therefore we say that $%
\Theta $ is invariant under the one-parameter point transformation with
generator, $X$, if \cite{Olver}
\begin{equation}
\mathcal{L}_{X^{\left[ 2\right] }}\Theta =0.  \label{go.11}
\end{equation}%
or equivalently%
\begin{equation}
\mathcal{L}_{X^{\left[ 2\right] }}\Theta =\lambda \Theta ~,~{mod}\Theta =0,
\label{go.12}
\end{equation}%
where $X^{\left[ 2\right] }$ is the second prolongation of $X$ in the space $%
J$. It is given by the formula
\begin{equation}
X^{\left[ 2\right] }=X+\eta _{i}\partial _{u_{,i}}+\eta _{ij}\partial
_{u_{,ij}},  \label{go.13}
\end{equation}%
where $\eta _{i}=D_{i}\left( \eta \right) -\Psi _{,k}D_{i}\left( \xi
^{k}\right) $, $\eta _{ij}=D_{j}\left( \eta _{i}\right) -u_{ki}D_{j}\xi ^{k}$
and $D_{i}$ is the operator of total differentiation, \textit{ie}, $D_{i}=%
\frac{\partial }{\partial x^{i}}+u_{,i}\frac{\partial }{\partial u}+u_{,ij}%
\frac{\partial }{\partial u_{,j}}+...~$\cite{Olver}. Moreover, if condition (%
\ref{go.11}) is satisfied (equivalently condition (\ref{go.12})), the vector
field $X$ is called a Lie point symmetry of the differential equation $%
\Theta $.

\subsection{Symmetries of linear evolution equations}

A geometric method which relates the Lie and the Noether point symmetries of
a class of second-order differential equations has been proposed in \cite%
{jgp,ijgmmp}. Specifically, the point symmetries of second-order partial
differential equations are related with the elements of the conformal
algebra of the underlying space which defines the Laplace operator.

Similarly, for the Lie symmetries of the second-order partial differential
equation,%
\begin{equation}
\Delta u-C^{\alpha }u_{,\beta }-u_{,t}=0,  \label{1fm.04}
\end{equation}%
where
\[
\Delta =\frac{1}{\sqrt{\left\vert g\right\vert }}\frac{\partial }{\partial
x^{\alpha }}\left( \sqrt{\left\vert g\right\vert }g^{\alpha \beta }\frac{%
\partial }{\partial x^{\beta }}\right)
\]
is the Laplace operator, $g_{\alpha \beta }=g\left( x^{\beta }\right) $ is a
nondegenerate tensor (we call it a metric tensor)\thinspace\ and $C^{\alpha
}=C^{\alpha }\left( x^{\beta }\right) $, the following theorem arises%
\footnote{%
For the proof see Appendix \ref{proof1}.}.

\begin{theorem}
\label{Theom1}The Lie point symmetries of (\ref{1fm.04}) are generated by
the Homothetic Group of the metric tensor $g_{\alpha \beta },~G_{H}$ which
defines the Laplace operator $\Delta $. The general form of the Lie symmetry
vector is
\begin{equation}
X_{L}=\left( c_{1}+2\psi _{I}\int T^{I}\left( t\right) dt\right) +\left(
T^{I}\left( t\right) Y_{I}^{\alpha }\left( x^{\beta }\right) \right)
\partial _{\alpha }+\left( a\left( x^{\beta },t\right) u+b\left( x^{\beta
},t\right) +c_{2}u\right) \partial _{u},  \label{1fm.05}
\end{equation}%
where $\psi _{I}$ is the homothetic factor of $Y_{I}^{a}$, $\psi _{I}=0$ for
the Killing vector (KV, $\psi _{I}=1$ for Homothetic vector (HV), $a\left(
x^{\beta },t\right) $ and $b\left( x^{\beta },t\right) $ are solutions of (%
\ref{1fm.05}), $Y_{I}^{\alpha }\left( x^{\beta }\right) $ is a KV/HV of $%
g_{\alpha \beta }$ and the following condition holds, namely,%
\begin{equation}
T^{I}L_{Y_{I}}C_{\alpha }-T_{,t}^{I}Y_{I\alpha }-2a_{,\alpha }=0.
\label{eq.18}
\end{equation}%
Note that $I=1,2,...,\dim G_{H}$.
\end{theorem}

Another important result for the $(1+1)$ linear evolution equation of the
form of (\ref{1fm.04}) is the following theorem which gives the dimension of
the possible admitted algebra.

\begin{theorem}
\label{Theom2} The one-dimensional linear evolution equation can admits 0,
1, 3 and 5 Lie point symmetries plus the homogenous and the infinity
symmetries\footnote{%
In the following the homogeneous and the infinity symmetries we call them
trivial symmetries.} \cite{lie}.
\end{theorem}

However, as equation (\ref{1fm.04}) is time independent, it admits always
the autonomous symmetry $\partial _{t}$. In the following we apply theorems %
\ref{Theom1} and \ref{Theom2} in order to study the Lie symmetries of the
space-dependent one-factor model.

\section{Space dependence of the one-factor model}

\label{space1}

The space-dependent one-factor model of commodity pricing is defined by the
equation%
\begin{equation}
\frac{1}{2}\sigma \left( S\right) ^{2}S^{2}F_{,SS}+\kappa \left( S\right)
\left( \mu \left( S\right) -\lambda \left( S\right) -\log S\right)
SF_{,S}-F_{,t}=0.  \label{1fm.06}
\end{equation}

The parameters, $\sigma$, $\kappa$, $\mu$ and $\lambda $, depend upon the
stock price, $S$. \ In order to simplify equation (\ref{1fm.06}) we perform
the coordinate transformation $S= \mbox{\rm e}^{x}$, that is, equation (\ref%
{1fm.06}) becomes%
\begin{equation}
\frac{1}{2}\sigma ^{2}\left( x\right) F_{,xx}+\left( \kappa \left( x\right)
\left( \mu \left( x\right) -\lambda \left( x\right) -x\right) -\frac{1}{2}%
\sigma ^{2}\left( x\right) \right) F_{,x}-F_{,t}=0  \label{1fm.07}
\end{equation}%
or
\begin{equation}
\Delta F+\left( \kappa \left( \mu -\lambda -x\right) -\frac{1}{2}\sigma ^{2}-%
\frac{1}{2}\sigma \sigma _{,x}\right) F_{,x}-F_{,t}=0,  \label{1fm.08}
\end{equation}%
where $\Delta $ is the Laplace operator in the one-dimensional space with
fundamental line element%
\begin{equation}
dx^{2}=\frac{2}{\sigma ^{2}\left( x\right) }dx^{2}  \label{1fm.09}
\end{equation}%
and admits a two-dimensional homothetic algebra. The gradient KV is $%
K^{1}=\sigma \left( x\right) \partial _{x}~$ and the gradient HV is $%
H^{2}=\sigma \left( x\right) \int \frac{1}{\sigma \left( x\right) }dx\ $with
homothetic factor$~\psi _{H}=1.$

Equation (\ref{1fm.08}) is of the form of (\ref{1fm.04}) where now
\begin{equation}
C^{x}\left( x\right) =-\left( \kappa \left( \mu -\lambda -x\right) -\frac{1}{%
2}\sigma ^{2}-\frac{1}{2}\sigma \sigma _{,x}\right)  \label{1fm.10}
\end{equation}%
and
\begin{equation}
C_{x}\left( x\right) =-\frac{2\left( \kappa \left( \mu -\lambda -x\right) -%
\frac{1}{2}\sigma ^{2}-\frac{1}{2}\sigma \sigma _{,x}\right) }{\sigma
^{2}\left( x\right) }.  \label{1fm10a}
\end{equation}%
Without performing any symmetry analysis we observe that, when $C^{x}\left(
x\right) =0$, (\ref{1fm.08}) is in the form of the heat conduction equation
and it is maximally symmetric, \textit{ie}, it admits $5+1+\infty $
symmetries. In the case for which $\sigma _{,x}\neq 0$ from (\ref{1fm.10})
we have that
\begin{equation}
\sigma ^{2}\left( x\right) =\left( 4\int e^{2x}K\left( x\right) \mbox{\rm d}
x+c_{1}\right) e^{-2x},  \label{1fm.11}
\end{equation}%
where $K\left( x\right) =\kappa \left( x\right) \left( \mu \left( x\right)
-\lambda \left( x\right) -x\right) $. However, this is only a particular
case whereas new cases can arise from the symmetry analysis.

\subsection{Symmetry analysis}

Let $Y_{I}\left( x\right) ,~I=1,2$ be the two HVs of the space (\ref{1fm.09}%
) with homothetic factors $\psi _{I}$. As (\ref{1fm.08}) is autonomous and
linear, it admits the Lie symmetries $\partial _{t},~F\partial _{F},~b\left(
t,x\right) \partial _{F}$, where $b\left( t,x\right) $ is a solution of (\ref%
{1fm.08}), Therefore from theorem \ref{Theom1} we have that the possible
additional Lie symmetry vector is
\begin{equation}
X=\left( 2\psi _{I}\int T^{I}\left( t\right) dt\right) +T^{I}\left( t\right)
Y_{I}\partial _{x}+\left( a\left( x,t\right) F\right) \partial _{F}
\end{equation}%
for which the following conditions hold%
\begin{equation}
T_{1}L_{Y_{1}}C_{x}+T_{2}L_{Y_{2}}C_{x}-T_{1,t}Y_{1}-T_{2,t}Y_{2}-2a_{,x}=0
\quad \mbox{\rm and}  \label{1fm.s04}
\end{equation}%
\begin{equation}
\Delta a-C^{x}a_{,x}-a_{,t}=0.  \label{1fm.s05}
\end{equation}

We study two cases: A) $a\left( t,x\right) =0$ and B) $a\left( t,x\right)
\neq 0.$

\subsubsection{Case A}

Let $a\left( t,x\right) =0$. Then (\ref{1fm.s05}) is satisfied. Hence from (%
\ref{1fm.s04}) we have the system
\begin{equation}
L_{Y_{I}}C_{x}-mY_{I}=0,  \label{1fm.s06}
\end{equation}%
where $T_{I,t}=mT$, \textit{ie}, $T_{I}\left( t\right) =T_{01}e^{mt}$. This
means that from any vector field $Y_{I}$ we have only one symmetry. Hence
from Theorem \ref{Theom2} condition (\ref{1fm.s06}) should hold for $I=1$
and $I=2.$ In this case the space-dependent one-factor model admits $%
3+1+\infty $ Lie point symmetries.

\subsubsection{Case B}

Consider that $a\left( t,x\right) \neq 0$. From (\ref{1fm.s04}) we have that%
\begin{equation}
a\left( t,x\right) =-\frac{1}{2}\int \left(
T_{1}L_{Y_{1}}C_{x}+T_{2}L_{Y_{2}}C_{x}-T_{1,t}Y_{1}-T_{2,t}Y_{2}\right)
dx+f\left( t\right)  \label{1fm.s061}
\end{equation}%
and then (\ref{1fm.s05}) gives%
\begin{eqnarray}
0 &=&\frac{1}{2}\sigma ^{2}\left(
T_{1}L_{Y_{1}}C_{x}+T_{2}L_{Y_{2}}C_{x}-T_{1,t}Y_{1}-T_{2,t}Y_{2,x}\right)
_{,x}+C^{x}\left(
T_{1}L_{Y_{1}}C_{x}+T_{2}L_{Y_{2}}C_{x}-T_{1,t}Y_{1}-T_{2,t}Y_{2}\right) +
\nonumber \\
&&-\int \left(
T_{1,t}L_{Y_{1}}C_{x}+T_{2,t}L_{Y_{2}}C_{x}-T_{1,tt}Y_{1}-T_{2,tt}Y_{2}%
\right) dx+2f\left( t\right) .  \label{1fm.s07}
\end{eqnarray}

Consider the case for which $T_{1}\neq T_{2}$. Recall that for the space (%
\ref{1fm.09}), $Y_{I,x}=2\psi _{I},$ that is, from (\ref{1fm.s07}) we have
the conditions%
\begin{eqnarray}
0 &=&\frac{1}{2}\sigma ^{2}\left( L_{Y_{1}}C_{x}\right)
_{,x}+C^{x}L_{Y_{1}}C_{x}-\frac{T_{1,t}}{T_{1}}\left( \frac{1}{2}\int
L_{Y_{1}}C_{x}dx+C^{x}Y_{1}+2\psi _{1}\right) +  \nonumber \\
&&+\frac{1}{2}\frac{T_{1,tt}}{T_{1}}\int Y_{1}dx+2f_{1,t}\quad
\mbox{\rm
and}  \label{1fm.s08}
\end{eqnarray}%
\begin{eqnarray}
0 &=&\frac{1}{2}\sigma ^{2}\left( L_{Y_{2}}C_{x}\right)
_{,x}+C^{x}L_{Y_{2}}C_{x}-\frac{T_{2,t}}{T_{2}}\left( \frac{1}{2}\int
L_{Y_{2}}C_{x}dx+C^{x}Y_{2}+2\psi _{2}\right) +  \nonumber \\
&&+\frac{1}{2}\frac{T_{2,tt}}{T_{2}}\int Y_{2}dx+2f_{2,t},  \label{1fm.s09}
\end{eqnarray}%
where $f\left( t\right) =f_{1}\left( t\right) +f_{2}\left( t\right) $. We
continue with the subcases:

\paragraph{Subcase B1}

Let $T_{I,tt}=T_{I}^{0}T_{I,t}~,~T_{I,t}\neq 0,$ that is, $T_{I}\left(
t\right) =T_{I}^{1}+T_{I}^{2}e^{T_{I}^{0}t}$. In this case the symmetry
conditions are:
\begin{equation}
0=\left( \frac{1}{2}\sigma ^{2}\left( L_{Y_{I}}C_{x}\right)
_{,x}+C^{x}L_{Y_{I}}C_{x}\right) -\frac{T_{I,t}}{T_{I}}\left( \frac{1}{2}%
\int L_{Y_{I}}C_{x}+C^{x}Y_{I}+2\psi _{I}+\frac{1}{2}\int Y_{I}\mbox{\rm d}%
x\right) +2f_{I,t}.  \label{1fm.s10}
\end{equation}%
Hence we have the following system
\begin{eqnarray}
\frac{1}{2}\sigma ^{2}\left( L_{Y_{I}}C_{x}\right) _{,x}+C^{x}L_{Y_{I}}C_{x}
&=&c  \label{1fm.s11} \\
\frac{1}{2}\int L_{Y_{I}}C_{x}+C^{x}Y_{I}+2\psi _{I}+\frac{1}{2}\int Y_{I}%
\mbox{\rm d}x &=&m  \label{1fm.s12} \\
2f_{I,t}-m\frac{T_{I,t}}{T_{I}}+c &=&0.  \label{1fm.s13}
\end{eqnarray}

If system (\ref{1fm.s11})-(\ref{1fm.s13}) holds for $I=1$ or $I=2$, then
equation (\ref{1fm.08}) admits $3+1+\infty $ $\ $ Lie symmetries and in the
case for which conditions (\ref{1fm.s11})-(\ref{1fm.s13}) hold, \textit{ie},
admits $5+1+\infty $ Lie symmetries which is the maximum for a $(1+1)$
evolution equation.

\paragraph{Subcase B2}

In the second subcase we consider that $T_{I,tt}\neq T_{I,t}$.

Hence, if B2.a) $T_{I,tt}=0$,~then from (\ref{1fm.s08}) and (\ref{1fm.s09})
it follows that%
\begin{equation}
\frac{1}{2}\int L_{Y_{I}}C_{x}dx+C^{x}Y_{I}+2\psi _{I}=0  \label{1fm.s14}
\end{equation}%
\begin{equation}
\frac{1}{2}\sigma ^{2}\left( L_{Y_{I}}C_{x}\right)
_{,x}+C^{x}L_{Y_{I}}C_{x}=c \quad \mbox{\rm and}  \label{1fm.s15}
\end{equation}%
\begin{equation}
2f_{I,t}+c=0,  \label{1fm.s16}
\end{equation}%
where from Theorem \ref{Theom2} these conditions must hold for $I=1$ and $%
I=2 $ and equation (\ref{1fm.08}) is maximally symmetric.

B2.b) Let~$T_{,tt}\neq 0$. Then\ it follows that%
\begin{equation}
\frac{1}{2}\int L_{Y_{I}}C_{x}dx+C^{x}Y_{I}+2\psi _{I}=0,  \label{1fm.s17}
\end{equation}%
\begin{equation}
\frac{1}{2}\sigma ^{2}\left( L_{Y_{I}}C_{x}\right) _{,x}+C^{x}L_{Y_{I}}C_{x}+%
\frac{1}{2}m\int Y_{I}dx=c \quad \mbox{\rm and}  \label{1fm.s18}
\end{equation}%
\begin{equation}
\frac{T_{I,tt}}{T_{I}}=m_{I}~,~2f_{I,c}+c-m_{I}=0.  \label{1fm.s19}
\end{equation}%
These conditions hold for $I=1$ or $I=2$. If these conditions hold for both $%
I=1$ and $2$, then equation (\ref{1fm.08}) is maximally symmetric.

We collect the results in the following theorem.

\begin{theorem}
\label{Theom3}The autonomous $(1+1)$ linear equation (\ref{1fm.08}), apart
from the symmetry of autonomy, the linear symmetry and the infinity
symmetry, can admit:

A) The two Lie symmetries $X_{I}=2\frac{\psi _{I}}{m}e^{mt}\partial
_{t}+e^{mt}Y_{I}\partial _{I}$, where $Y_{I}$ is a HV of the one-dimensional
flat space with $I=`1,2$ if and only if condition (\ref{1fm.s06}) holds for $%
I=1$ and $I=2$.

B1) The two or four Lie symmetries
\begin{equation}
X_{I}=2\psi _{I}\int T_{I}\left( t\right) dt~\partial _{t}+T_{I}\left(
t\right) Y_{I}+\alpha \left( t,x\right) F\partial _{F}  \label{1fm.s20a}
\end{equation}
if conditions (\ref{1fm.s11})-(\ref{1fm.s13}) hold for $I=1$ or $2,$ and $%
I=1 $ and $2$, respectively, where $T_{I}\left( t\right)
=T_{I}^{1}+T_{I}^{2}e^{T_{I}^{0}t}$ and
\begin{equation}
a\left( t,x\right) =-\frac{1}{2}\int \left(
T_{I}L_{Y_{I}}C_{x}-T_{I,t}Y_{1}\right) dx+f_{I}\left( t\right).
\label{1fm.s20}
\end{equation}

B2.a) The four Lie symmetries (\ref{1fm.s20a}) if conditions (\ref{1fm.s14}%
)--(\ref{1fm.s16}) hold for $I=1$ and $2$, where $T_{I}=T_{I0}+T_{I1}t$ and $%
a\left( t,x\right) $ is given by (\ref{1fm.s20}).

B2.b) The two or four Lie symmetries (\ref{1fm.s20a}) if and only if
conditions (\ref{1fm.s17})--(\ref{1fm.s19}) hold for $I=1$ or $2,$ and $I=1$
and $2$, respectively, where $T_{I,tt}=m_{I}T.$
\end{theorem}

Furthermore, we comment that theorem (\ref{Theom3}) holds for all linear
autonomous equations of the form of (\ref{1fm.04}).

Here we discuss the relation among the Lie symmetries and the initial
condition (\ref{1fm.01a}). In the case of constant parameters, \textit{ie},
in equation (\ref{1fm.01}) the Lie symmetry vector (\ref{1fm.03}) is the
linear combination among the linear symmetry $F\partial _{F} $ and the
symmetry which is generated by the KV of the underlying space, which is $%
K^{1}=\sigma _{0}\partial _{x},~$for $\sigma \left( x\right) =\sigma _{0}$.
However, for a general function, $\sigma \left( S\right) $, in order for the
symmetry which is generated by the KV $K^{1}$ to satisfy the initial
condition $\sigma \left( S\right) =\sigma _{0}$ or the initial condition has
to change.

Consider now that $Y^{1}=K^{1}$ and satisfies the conditions
\begin{equation}
L_{K_{1}}C_{x}=0~,~C^{x}Y_{1}=0.
\end{equation}%
Then from theorem \ref{Theom3}, B2.a, we have that $\sigma \left( x\right) $
is given by (\ref{1fm.11}) and at the same time $Y^{2}=H$ generates two Lie
point symmetries for equation (\ref{1fm.08}). The Lie point symmetries are
\[
X_{t}=\partial _{t},~X_{F}=F\partial _{F},~X_{1}=\sigma \left( x\right)
\partial _{x},~X_{2}=t\sigma \left( x\right) \partial _{x}-\sigma \left(
x\right) F\partial _{F},~
\]%
\[
X_{3}=2t\partial _{t}+\sigma \left( x\right) \int \frac{1}{\sigma \left(
x\right) }dx~\partial _{x} \mbox{\rm and}
\]%
\[
X_{4}=t^{2}\partial _{t}+t\sigma \left( x\right) \int \frac{1}{\sigma \left(
x\right) }dx~\partial _{x}-\left( \frac{1}{2}t+\left( \int \frac{1}{\sigma
\left( x\right) }dx\right) ^{2}\right) F\partial _{F}
\]%
plus the autonomous and trivial symmetries. The symmetry vector field $X_{1}$
is the KV of the one-dimensional space. Therefore, if we wish the field $%
\bar{X}=X_{1}+\mu X_{F}$ to satisfy an initial condition such as $F\left(
0,x\right) =g\left( x\right) $, then it should be $X_{1}\left( g\left(
x\right) \right) =g\left( x\right) $ which gives $g\left( x\right) =e^{\int
\sigma \left( x\right) dx}$. \ From this we can see that, when $\sigma
\left( x\right) =\sigma _{0}$, we have the initial condition (\ref{1fm.01a}).

Let $\kappa $, $\lambda $ and $\mu $ be constants. Hence from (\ref{1fm.11})
we have that
\begin{equation}
\sigma ^{2}\left( x\right) =2\kappa \left( \left( \mu -\lambda \right)
-x\right) +\kappa +c_{1}e^{-2x},
\end{equation}%
where for $c_{1}=0$ we have
\begin{equation}
\ln g\left( x\right) =\mp \kappa \frac{\left( 2\left( \mu -\lambda -x\right)
+1\right) }{3}\sigma
\end{equation}%
and the solution for position $\sigma \left( x\right) >0$%
\begin{equation}
\ln F\left( t,x\right) =\frac{\mu }{2}\left( \mu \kappa t-2\sqrt{2\left( \mu
-\lambda -x\right) +1}\right) .
\end{equation}

Let now $\sigma \left( x\right) =x$ and consider that the KV $%
K^{1}=x\partial _{x}$ generates a Lie point symmetry of equation (\ref%
{1fm.08}) from Case A of theorem \ref{Theom3}. Then from condition (\ref%
{1fm.s06}) we have that%
\begin{equation}
xC_{,x}^{x}-C^{x}+mx=0,
\end{equation}%
that is,
\begin{equation}
C^{x}=-m\ln x+c_{1}x.
\end{equation}

However, in that case, equation (\ref{1fm.08}) is maximally symmetric and
admits $5+1+\infty ~$Lie point symmetries. Consider reduction with the Lie
symmetry $\bar{X}=e^{mt}x\partial _{x}+F\partial _{F}$ which keeps invariant
the initial condition%
\begin{equation}
F\left( 0,x\right) =\frac{1}{x}.
\end{equation}

The application of $\bar{X}$ in (\ref{1fm.08}) gives
\begin{equation}
F\left( t,x\right) =x^{\exp \left( -mt\right) }\exp \left( -\frac{1}{4m}%
e^{-2mt}-\frac{c_{1}}{m}e^{-mt}\right) ~,~m\neq 0
\end{equation}%
\begin{equation}
F\left( t,x\right) =x\exp \left( \left( c_{1}+\frac{1}{2}\right) t\right)
~,~m=0.
\end{equation}

As another application of theorem \ref{Theom3} we select $\sigma \left(
x\right) =e^{x}$. Then the KV $K^{1}$ is $K^{1}=e^{x}\partial _{x}$. Let
this generate a Lie point symmetry for equation\ (\ref{1fm.08}) from the
case A of theorem \ref{Theom3}, that is, conditions (\ref{1fm.s06}) give%
\begin{equation}
C^{x}=-m+c_{1}e^{-x},
\end{equation}%
where now we can see that equation (\ref{1fm.08}) is maximally symmetric and
admits $5+1+\infty $ point symmetries. Consider the Lie symmetry $\bar{X}%
=e^{mt+x}\partial _{x}+F\partial _{F}$, which leaves invariant the modified
initial condition $F\left( 0,x\right) =x$. The invariant solution which
follows is
\[
\ln \left( F,x\right) =-\frac{e^{-mt}}{4m}\left( 4me^{-x}-4c+e^{-mt}\right)
.
\]

We observe that, when $K^{1}$ generates a Lie point symmetry for equation (%
\ref{1fm.08}), the functional form of $C^{x}$, which includes $\kappa \left(
x\right) ,\lambda \left( x\right) $ and $\mu \left( x\right) $ has a
specific form, such that equation (\ref{1fm.08}) is maximally symmetric and
equivalent with the Black-Scholes and the Heat equations. In general, for
unknown function $\sigma \left( x\right) $, from theorem \ref{Theom1} we
have the following corollary.

\begin{corollary}
\label{Cor}When the KV of the underlying space which defines the Laplace
operator in equation (\ref{1fm.08}) generates a Lie point symmetry, the
functional form of $C^{x}$ is
\begin{equation}
C^{x}\left( x\right) =\sigma \left( x\right) \int \frac{m}{\sigma \left(
x\right) }dx+c\sigma \left( x\right)
\end{equation}%
and equation (\ref{1fm.08}) is maximally symmetric. The symmetry vectors,
among the autonomous, the homogeneous and the infinity symmetries, are:%
\begin{equation}
Z^{1}=e^{mt}K^{1}~,~Z^{2}=e^{-mt}\left( K^{1}+m\int \frac{dx}{\sigma \left(
x\right) }F\partial _{F}\right)
\end{equation}%
\begin{equation}
Z^{3}=e^{2mt}\left( \partial _{t}+H\right)
\end{equation}%
\begin{equation}
Z^{4}=e^{-2mt}\left( \partial _{t}-H-m\left( 2m\left( \int \frac{dx}{\sigma
\left( x\right) }\right) ^{2}-1\right) F\partial _{F}\right)
\end{equation}%
for $m\neq 0$, $c=0$,
\begin{equation}
Z^{1}=e^{mt}K_{1}~,~Z^{2}=e^{-mt}\left( K^{1}+\left( m\int \frac{dx}{\sigma
\left( x\right) }+c\right) F\partial _{F}\right)
\end{equation}%
\begin{equation}
Z^{3}=e^{2mt}\left( \partial _{t}+H+cK^{1}\right)
\end{equation}%
\begin{equation}
Z^{4}=e^{-2mt}\left( \partial _{t}-H-cK^{1}+\left( 2m^{2}\left( \int \frac{dx%
}{\sigma \left( x\right) }\right) ^{2}+4mc\left( \int \frac{dx}{\sigma
\left( x\right) }\right) +2c^{2}-m\right) F\partial _{F}\right)
\end{equation}%
for $m\neq 0$,~$c\neq 0,$ and%
\begin{equation}
Z^{1}=K^{1}~,~Z^{2}=tK^{1}-\left( \int \frac{dx}{\sigma \left( x\right) }%
F-t\right) F\partial _{F}
\end{equation}%
\begin{equation}
Z^{3}=2t\partial _{t}+H+ctK^{1}
\end{equation}%
\begin{equation}
Z^{4}=t^{2}\partial _{t}+H+\frac{t}{2c}K^{1}-\left( \frac{1}{2}\int \frac{dx%
}{\sigma \left( x\right) }^{2}+\frac{1}{2c}\int \frac{dx}{\sigma \left(
x\right) }-ct\int \frac{dx}{\sigma \left( x\right) }+\frac{c^{2}}{2}%
t^{2}\right) F\partial _{F}
\end{equation}%
for $m=0$, $c\neq 0$, where $K_{1}=\sigma \left( x\right) \partial _{x}$ and
$H=\sigma \left( x\right) \int \frac{dx}{\sigma \left( x\right) }\partial
_{x}$ are the elements of the Homothetic algebra of the underlying space.
\end{corollary}

We note that corollary \ref{Cor} holds for all autonomous linear 1+1
evolution equations. In the following section we discuss the group
invariants of the time-dependent problem.

\section{Time-dependent one-factor model}

\label{proof2}

When the parameters $\sigma ,~\kappa ,\lambda $ of equation (\ref{1fm.01})
depend upon time, the one-factor model can be written as
\begin{equation}
\frac{1}{2}\sigma ^{2}(t)(F_{xx}-F_{x})+(p(t)-xq(t))F_{x}-F_{t}=0,
\label{A2}
\end{equation}%
where%
\[
p\left( t\right) =q\left( t\right) \left( \mu \left( t\right) -\lambda
\left( t\right) \right) ~,~q\left( t\right) =\kappa \left( t\right) .
\]

Without loss of generality we can select $\sigma \left( t\right) =1$. \ By
analysing the determining equations as provided by the Sym package \cite%
{Dimas05a,Dimas06a,Dimas08a} we find that the general form of the Lie
symmetry vector is
\begin{eqnarray}
X &=&a(t)\partial _{t}+\left[ b(t)+\frac{a^{\prime }x}{2}\right] \partial
_{x}  \nonumber  \label{A3} \\[1pt]
&&+\left[ f(t)+\frac{1}{4}\left( 4xbq+x(1-2p)a^{\prime 2}qa^{\prime
}-4x(b^{\prime }+ap^{\prime })\right. \right.   \nonumber \\[1pt]
&&\left. \left. -x^{2}(-2aq^{\prime }+a^{\prime \prime })\right) \right]
F\partial _{F},
\end{eqnarray}%
where functions $a\left( t\right) ,b\left( t\right) \,,~f\left( t\right) $
are given by the system of ordinary differential equations,
\begin{eqnarray}
0 &=&-bq+2bpq-\frac{a^{\prime }}{4}+pa^{\prime 2}a^{\prime }+qa^{\prime }+
\nonumber  \label{A4} \\[1pt]
&&+b^{\prime }-2pb^{\prime }-2f^{\prime }+ap^{\prime }-2app^{\prime
}+aq^{\prime }-\frac{a^{\prime \prime }}{2},
\end{eqnarray}%
\begin{eqnarray}
0 &=&-2bq^{2}-\frac{3}{2}qa^{\prime }+3pqa^{\prime }+2aqp^{\prime }+
\nonumber  \label{A5} \\[1pt]
&&+3a^{\prime }p^{\prime }-aq^{\prime }-2bq^{\prime }+2apq^{\prime
}+2b^{\prime \prime }+2ap^{\prime \prime }
\end{eqnarray}%
and
\[
0=-2q^{2}a^{\prime }-2aqq^{\prime }-2a^{\prime }q^{\prime }-aq^{\prime
\prime }+\frac{a^{\prime \prime \prime }}{2}
\]%
in addition to the infinite number of solution symmetries. Consequently the
algebra is $\{sl(2,\Re )\oplus _{s}W_{3}\}\oplus _{s}\infty A_{1}$ so that
it is related to the classical Heat Equation by means of a point
transformation. In the following we discuss our results.

\section{Conclusions}

\label{con}

In the models of financial mathematics the parameters of the models are
assumed to be constants. However, in real problems these parameters can
depend upon the stock prices and upon time. In this work we considered the
one-factor model of Schwartz and we studied the Lie symmetries in the case
for which the parameters of the problem are space-dependent. In terms of Lie
symmetries, the one-factor model it is maximally symmetric and it is
equivalent with the Heat equation, but in the case where the parameters are
space dependent, that is not necessary true, and we show that the model can
admit 1, 3 or 5 Lie point symmetries (except the trivial ones). To perform
this analysis we studied the Lie symmetries of the autonomous linear
evolution equation and we found that there exist a unique relation among the
Lie symmetries and the collineations of the underlying geometry, where as
geometry we define the \textquotedblleft space\textquotedblright\ of the
second derivatives. However, for a specific relation among the parameters of
the model the system is always maximally symmetric. \ In particular, that
holds when $\sigma \left( x\right) $ is an arbitrary function and
\begin{equation}
\left( \kappa \left( \mu -\lambda -x\right) -\frac{1}{2}\sigma ^{2}-\frac{1}{%
2}\sigma \sigma _{,x}\right) =\sigma \int \frac{m}{\sigma }dx+c\sigma,
\label{con.01}
\end{equation}%
where $m,c$ are constants. In that case, the correspoding symmetry (\ref%
{1fm.03}) becomes $\bar{Z}=e^{mt}\sigma \left( x\right) \partial
_{x}+F\partial _{x}$.

Consider that $\sigma \left( x\right) =1+\varepsilon e^{x-x_{0}}$, and (\ref%
{con.01}) holds. Then the application of the Lie symmetry $\bar{Z}$ in (\ref%
{1fm.08}) gives the solution%
\begin{equation}
\ln F\left( t,x\right) =-e^{-mt}\ln \left( 1+\varepsilon e^{x-x_{0}}\right)
+e^{-mt}\left( x-x_{0}\right) +\frac{c}{m}e^{-mt}-\frac{1}{4m}e^{-2mt},
\label{con.02}
\end{equation}%
where in the limit $\varepsilon \rightarrow 0$, solution (\ref{con.02})
becomes%
\begin{equation}
\ln F\left( t,x\right) =e^{-mt}\left( x-x_{0}\right) +\frac{c}{m}e^{-mt}-%
\frac{1}{4m}e^{-2mt}  \label{con.03}
\end{equation}%
which can compared with solution (\ref{1fm.02}).

Consider now that $\sigma \left( x\right) $ is periodic around the line $%
\sigma _{0}=1$. Let $\sigma \left( x\right) =1+\varepsilon \sin \left(
\omega x\right) $ that and (\ref{con.01}) holds. Hence the solution of the
space-dependent one-factor model (\ref{1fm.08}) which follows from the Lie
symmetry $\bar{Z}$ is
\begin{equation}
\ln F\left( t,x\right) =\frac{2e^{-mt}}{\sqrt{1+\varepsilon ^{2}}}i\arctan
\left( \frac{\tan \left( \frac{x}{2}\right) +\varepsilon }{\sqrt{%
1-\varepsilon ^{2}}}\right) +\frac{c}{m}e^{-mt}-\frac{1}{4m}e^{-2mt}
\label{con.04}
\end{equation}%
which is a periodic function of the stock price $x$. \ For $\varepsilon <<1$
the Taylor expansion of the static solution (\ref{con.04}) around the point $%
\varepsilon =0$, is%
\begin{equation}
\ln F\left( t_{0},x\right) \simeq ix+\varepsilon \cos x+O\left( \varepsilon
^{2}\right) .
\end{equation}

\begin{figure}[tbp]
\centering
\includegraphics[height=7cm]{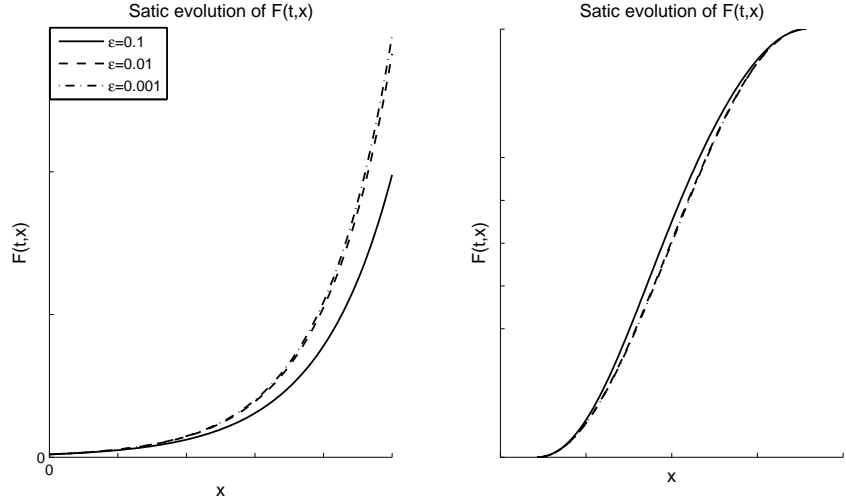}
\caption{Evolution of the static solution (\protect\ref{con.02}) (left
figures) and solution (\protect\ref{con.04}) (right. figure) for various
values of the constant $\protect\varepsilon $. Solid line is for $\protect%
\varepsilon =0.1$, dash dash line is for $\protect\varepsilon =0.1$, and the
dash dot line is for $\protect\varepsilon =0.01$.}
\label{fig3}
\end{figure}

In Figure \ref{fig3} we give the static evolution of the solutions, (\ref%
{con.03}) and (\ref{con.04}), for various values of the constant $%
\varepsilon $.

On the other hand, in Section \ref{proof2} we studied the case for which the
parameters of the one-factor model are time-dependent and we showed that the
model is always maximally symmetric and equivalent with the Heat Equation,
that is, the time-depedence does not change the admitted group invariants of
the one-factor model (\ref{1fm.01}).

A more general consideration will be to extend this analysis to the
two-factor and three-factor models and also to study the cases for which the
parameters are dependent upon the stock price and upon the time, \textit{ie}%
, the parameters are space and time dependent. This work is in progress.

Finally we remark how useful are the methods which are applied in Physics
and especially in General Relativity for the study of space-dependent
problems in financial mathematics. The reason for this is that from the
second derivatives a (pseudo)Riemannian manifold can be defined. This makes
the use of the methods of General Relativity and Differential Geometry
essential.

\subsubsection*{Acknowledgements}

The
research of AP was supported by FONDECYT postdoctoral grant no. 3160121. RMM thanks the National Research Foundation of the Republic of
South Africa for the granting of a postdoctoral fellowship with grant number
93183 while this work was being undertaken.

\appendix{}

\section{Proof of Theorem \protect\ref{Theom1}}

\label{proof1}

In \cite{jgp} it has been shown that for a second-order PDE of the form,
\begin{equation}
A^{ij}\left( x^{k}\right) u_{,ij}-B^{k}(x)u_{k}=0,  \label{eq.02}
\end{equation}%
the Lie symmetries are generated by the conformal algebra of the tensor $%
A^{ij}$. \ Specifically the Lie symmetry conditions for equation (\ref{eq.02}%
) are
\begin{equation}
A^{ij}(a_{ij}u+b_{ij})-(a_{,i}u+b_{,i})B^{i}=0,  \label{eq.03}
\end{equation}%
\begin{equation}
A^{ij}\xi _{,ij}^{k}-2A^{ik}a_{,i}+aB^{k}-\xi _{,i}^{k}B^{i}+\xi
^{i}B_{,i}^{k}-\lambda B^{k}=0\quad \mbox{\rm and}  \label{eq.04}
\end{equation}%
\begin{equation}
L_{\xi ^{i}}A^{ij}=(\lambda -a)A^{ij},  \label{eq.05}
\end{equation}%
where%
\begin{equation}
\eta =a(x^{i})u+b(x^{i})~~,~\xi ^{k}=\xi ^{k}(x^{i}).  \label{eq.07}
\end{equation}

By comparison of equations (\ref{1fm.04}) and (\ref{eq.02}) we have that~$%
x^{i}=\left( x^{\alpha },t\right) ,~A^{ij}=g^{\alpha \beta }\left( x^{\beta
}\right) $, \textit{ie}, $A^{tt}=A^{it}=0~$ and%
\begin{equation}
B^{t}=1~,~B^{\alpha }=\Gamma ^{\alpha }+C^{\alpha },
\end{equation}%
where $\Gamma ^{\alpha }=\Gamma _{\beta \gamma }^{\alpha }g^{\beta \gamma }$%
. Therefore the symmetry vector $X$ for equation (\ref{1fm.04}) has the form
\begin{equation}
X=\xi ^{t}\left( t\right) \partial _{t}+\xi ^{\alpha }\left( t,x^{\beta
}\right) \partial _{\alpha }+\left( a(x^{\alpha },t)u+b(x^{a},t)\right)
\partial _{u}.  \label{eq.08}
\end{equation}%
We continue with the solution of the symmetry conditions.

When we replace $~A^{ij}=g^{\alpha \beta }\left( x^{\beta }\right) $~ in (%
\ref{eq.05}), it follows that $L_{\xi ^{\alpha }}g_{\alpha \beta }=\left(
\alpha -\lambda \right) g_{\alpha \beta }~$ which means that $\xi ^{\alpha
}=T^{I}\left( t\right) Y_{I}^{\alpha }\left( x^{\beta }\right) $, where $%
Y_{I}^{\alpha }$ is a CKV of the metric, $g_{\alpha \beta }$, with conformal
factor $\psi _{I}$,$~$ \textit{ie},$~\psi _{I}=\frac{1}{n}Y_{;\alpha
}^{\alpha }$ and $\alpha -\lambda =2T^{I}\psi _{I}.~$ Furthermore, from (\ref%
{eq.04}) the following system follows (recall that $B_{,i}^{t}=0$~and $\xi
_{,\beta }^{t}=0\,$)
\begin{equation}
g^{ij}\xi _{,ij}^{\alpha }-2g^{i\alpha }a_{,i}+aB^{\alpha }-\xi
_{,i}^{\alpha }B^{i}+\xi ^{i}B_{,i}^{\alpha }-\lambda B^{\alpha }=0\quad %
\mbox{\rm and}  \label{eq.11}
\end{equation}%
\begin{equation}
\left( a-\lambda \right) B^{t}-\xi _{,t}^{t}B^{t}=0.  \label{eq.12}
\end{equation}%
Moreover we observe that~$\psi \left( x^{k}\right) =\psi _{I}$, where $\psi
_{I}$ is a constant; that is, $Y_{I}^{\alpha }\left( x^{\beta }\right) $ is
a KV/HV of $g_{\alpha \beta }$. \ Finally for the function, $\xi ^{t}\left(
t\right) ,$ it holds that $\xi ^{t}\left( t\right) =2\psi _{I}\int
T^{I}\left( t\right) dt.~$

The symmetry condition (\ref{eq.11}) gives
\begin{equation}
g^{\beta \gamma }\xi _{,\beta \gamma }^{\alpha }-\xi _{,\gamma }^{\alpha
}\Gamma ^{\gamma }+\xi ^{\gamma }\Gamma _{,\gamma }^{\alpha }-\xi _{,\gamma
}^{\alpha }C^{\gamma }+\xi ^{\gamma }C_{,\gamma }^{\alpha }-2g^{\beta \alpha
}a_{,\beta }+\left( a-\lambda \right) \Gamma ^{\alpha }+\left( a-\lambda
\right) C^{\alpha }-\xi _{,t}^{\alpha }=0.  \label{eq.15}
\end{equation}

It is well known that
\begin{equation}
g^{\beta \gamma }\xi _{,\beta \gamma }^{\alpha }+\xi ^{\gamma }\Gamma
_{,\gamma }^{\alpha }-\xi _{,\gamma }^{\alpha }\Gamma ^{\gamma }+\left(
a-\lambda \right) \Gamma ^{\alpha }=g^{\beta \gamma }\left( L_{\xi ^{\alpha
}}\Gamma _{\beta \gamma }^{\alpha }\right)  \label{eq.16}
\end{equation}%
and, as $\xi ^{\alpha }$ is a HV of $g_{\alpha \beta }$, $\left( L_{\xi
^{\alpha }}\Gamma _{\beta \gamma }^{\alpha }\right) =0$ holds. Therefore (%
\ref{eq.15}) becomes%
\begin{equation}
T^{I}\left( L_{Y_{I}^{\alpha }}C^{\alpha }+2\psi _{I}C^{\alpha }\right)
-T_{,t}^{I}Y_{I}^{\alpha }-2g^{\beta \alpha }a_{,\beta }=0.  \label{eq.17}
\end{equation}

However, because $g_{\alpha \gamma }L_{\xi }C^{\alpha }=\delta _{\gamma
}^{\beta }L_{\xi }C_{\beta }-2\psi C_{\gamma },$ condition (\ref{eq.17}) can
be written as%
\begin{equation}
T^{I}L_{Y_{I}}C_{a}-T_{,t}^{I}Y_{I\alpha }-2a_{,\beta }=0.  \label{eq.18a}
\end{equation}

Finally from condition (\ref{eq.03}) we have the system
\begin{eqnarray}
\Delta a-C^{\alpha }a_{,a}-a_{,t} &=&0 \quad \mbox{\rm and}  \label{eq.19} \\
\Delta b_{\alpha \beta }-C^{\alpha }b_{,a}-b_{,t} &=&0  \label{eq.20a}
\end{eqnarray}%
which means that $a\left( x^{k},t\right) $ and $b\left( x^{k},t\right) $ are
solutions of (\ref{1fm.04}). We continue with the study of some special
cases:

Case I: Let $T^{I}\left( t\right) =0$. Then from (\ref{eq.18}) $a_{,\beta
}=0 $ which means that $\alpha =a\left( t\right) $. However, from (\ref%
{eq.19}) we have that $\alpha \left( t\right) =a_{0}$ which gives the linear
symmetry $a_{0}u\partial _{u}$. In that case from the form of $\xi ^{t}$ the
autonomous symmetry $\partial _{t}~$ arises.

Case II: For $T^{I}\left( t\right) \neq 0~$\ the generic symmetry vector is%
\begin{equation}
X=2\psi _{I}\int T^{I}\left( t\right) \partial _{t}+T^{I}Y_{I}^{\alpha
}\partial _{\alpha }+\left( a\left( x^{\beta },t\right) \right) u\partial
_{u},
\end{equation}%
where conditions (\ref{eq.18a}) and (\ref{eq.19}) hold.

\end{document}